\documentclass{amsart}
\usepackage{amsmath}
\usepackage{amsfonts}
\usepackage{amssymb}
\usepackage{graphicx}%
\usepackage{mathrsfs}

\newtheorem{theorem}{Theorem}

\newtheorem*{GGR+theorem}{GGR+ Theorem}

\newtheorem{conjecture}[theorem]{Conjecture}
\newtheorem{corollary}[theorem]{Corollary}

\newtheorem{example}[theorem]{Example}
\newtheorem{exercise}[theorem]{Exercise}
\newtheorem{lemma}[theorem]{Lemma}

\newtheorem{proposition}[theorem]{Proposition}

\numberwithin{theorem}{section}
\numberwithin{equation}{section}

\newcommand{\ot}{\otimes}

\def\J{{\mathbb J}}
\def\I{{\mathbb I}}

\def\Z{{\mathbb Z}}

\newcommand{\cD}{\mathcal{D}}
\newcommand{\ttf}{\mathtt{f}}
\newcommand{\ttk}{\mathtt{k}}
\newcommand{\itemo}{\item[{}]}
\def\Hom {{\operatorname{Hom}}}

\def\N{{\mathbb N}}

\def\H{{\mathscr H}}

\newcommand{\itema}{\item[{{\rm(a)}}]}
\newcommand{\itemb}{\item[{{\rm(b)}}]}
\newcommand{\itemc}{\item[{{\rm(c)}}]}
\newcommand{\itemd}{\item[{{\rm(d)}}]}
\newcommand{\iteme}{\item[{{\rm(e)}}]}

\newcommand{\cS}{\mathcal{S}}

\newcommand{\cH}{\mathcal{H}}

\newcommand{\ff}{\footnote}
\newcommand{\noi}{\noindent}
\newcommand{\ga}{\alpha}
\newcommand{\gb}{\beta}
\newcommand{\gc}{\gamma}

\newcommand{\Gd}{\Delta}
\newcommand{\Gt}{\Theta}
\newcommand{\gd}{\delta}

\newcommand{\gs}{\sigma}

\newcommand{\gl}{\lambda}

\newcommand{\gr}{\rho}
\newcommand{\gk}{\kappa}
\newcommand{\gep}{\epsilon}
\newcommand{\gth}{\theta}

\newcommand{\A}{\mathbb A}
\newcommand{\lra}{\longrightarrow}
\newcommand{\fgl}{\mathfrak{gl}}
\newcommand{\fn}{\mathfrak{n}}
\newcommand{\fb}{\mathfrak{b}}

\newcommand{\fh}{\mathfrak{h}}
\newcommand{\ey}{\end{eqnarray}}
\newcommand{\by}{\begin{eqnarray}}
\newcommand{\nn}{\nonumber}

\newcommand{\bco}{\begin{conjecture}}
\newcommand{\ba}{\begin{alg}}
\newcommand{\ea}{\end{alg}}
\newcommand{\eco}{\end{conjecture}}
\newcommand{\bpf}{\begin{proof}}
\newcommand{\epf}{\end{proof}}
\newcommand{\bt}{\begin{theorem}}

\newcommand{\et}{\end{theorem}}

\newcommand{\br}{\begin{rem}}
\newcommand{\er}{\end{rem}}
\newcommand{\brs}{\begin{rems}}
\newcommand{\ers}{\end{rems}}
\newcommand{\bi}{\begin{itemize}}
\newcommand{\ei}{\end{itemize}}
\newcommand{\bl}{\begin{lemma}}
\newcommand{\bsul}{\begin{sublemma}}
\newcommand{\esul}{\end{sublemma}}
\newcommand{\bp}{\begin{proposition}}
\newcommand{\be}{\begin{equation}}
\newcommand{\bc}{\begin{corollary}}
\newcommand{\bexs}{\begin{examples}}
\newcommand{\eexs}{\end{examples}}
\newcommand{\bexa}{\begin{example}}
\newcommand{\eexa}{\end{example}}
\newcommand{\bex}{\begin{exercise}}
\newcommand{\eex}{\end{exercise}}
\newcommand{\btab}{\begin{tab}}
\newcommand{\etab}{\end{tab}}
\newcommand{\bg}{\begin{fig}}

\newcommand{\gL}{\Lambda}

\newcommand{\fm}{\mathfrak{m}}

\newcommand{\fk}{\mathfrak{k}}

\newcommand{\sfg}{\small{\mathfrak{g}}}

\def\ad{{\operatorname{ad}}}
\newcommand{\eg}{\end{fig}}
\newcommand{\el}{\end{lemma}}
\newcommand{\ep}{\end{proposition}}
\newcommand{\ee}{\end{equation}}
\newcommand{\ec}{\end{corollary}}
\newcommand{\Bc}{\begin{center}}
\newcommand{\Ec}{\end{center}}
\newcommand{\bh}{\begin{hyp}}
\newcommand{\eh}{\end{hyp}}
\newcommand{\bhs}{\begin{hyps}}
\newcommand{\ehs}{\end{hyps}}
\newcommand{\bd}{\begin{dfn}}
\newcommand{\ed}{\end{dfn}}
  \newtheorem{rems}{Remarks}
  \newtheorem{rem}{Remark}

\numberwithin{equation}{section}%
\begin{document}
\newcommand{\fg}{\mathfrak{g}}
\newcommand{\fsl}{\mathfrak{sl}}

\title{Shapovalov Elements for $U_q(\fsl(N+1))$}
\author{Stefan Catoiu}
\address{Department of Mathematics, DePaul University, Chicago, IL 60614}
\email{scatoiu@depaul.edu}

\author{Ian M. Musson}
\address{Department of Mathematics, University of Wisconsin-Milwaukee, Milwaukee, WI }
\email{musson@uwm.edu}


\begin{abstract}
For a simple Lie algebra, Shapovalov elements give rise to highest weight vectors in Verma modules. The usual construction of these elements uses induction on the length of a certain Weyl group element.  If $\fg= \fsl(N+1)$ explicit expressions for Shapovalov elements were given in \cite{M2}. Here we adapt the argument to the quantized enveloping algebra of $\fg$.

\end{abstract}
\maketitle
\section{The Quantized Enveloping Algebra}
\subsection{The PBW Theorem}
Let $\ttf$ be a field and $\ttk = \ttf(q)$ the field of rational functions over $\ttf$.  Unadorned tensors are taken over $\ttk.$ Let $\fg= \fsl(N+1)$ and   
 $A = (a_{i,j})$ the Cartan matrix of $\fg$. Let $U_q(\fg) $ be the $\ttk$-algebra with generators $e_i, f_i,  k_i^{\pm 1}$ for $i\in[N]:=\{1,\ldots ,N\}$. As the notation suggests,  $k_i$ is a unit in  $U_q(\fg) $ with two-sided inverse $k_i^{-1}$. Furthermore  $U_q(\fg) $ satisfies the following relations for all $i$ and $j$:
\be \label{yak} k_ik_j,  =k_jk_i,\ee
\be \label{eq1.1} k_ie_jk_i^{-1}=q^{a_{i,j}}e_j,  \;\;k_if_jk_i^{-1}=q^{-a_{i,j}}f_j, \ee
\begin{equation}\label{owl}
e_if_j-f_je_i=\delta_{ij}(k_i^2-k_i^{-2})/(q^2-q^{-2}).
\end{equation}
as well as the Serre relations, \cite{Y}, (1.4) and (1.5) which we do not use. They are needed for the proof of Yamane's PBW theorem, Theorem \ref{YT} (c).  We denote the subalgebra of $U_q(\fg)$ generated by  $ f_i$  (resp.  $ e_i$) for $i\in[N]$ by 
$U_q(\fn^-)$ (resp. $U_q(\fn^+)$).\\ \\
Define $f_{i,i+1}= f_i$.  This is a quantum analog of the negative root vector $E_{i+1,i}\in \fsl(N+1)$. For analogs of other negative root vectors,  
consider the following elements, introduced by Jimbo in \cite{Ji}. For $j>i+1$, define
\begin{equation}\label{dog}
f_{i,j}=qf_{i,j-1}f_{j-1,j}-q^{-1}f_{j-1,j}f_{i,j-1}.
\end{equation}
Define the lexicographic order $<$ on 
 $\N \times \N$ by
$(i,j) < (m, n)$ iff $i < m$ or  $i = m, j < n.$ \bt \label{YT} 
\bi \itema The group of units $\A$ of $U_q(\fg)$ generated by the $ k_i^{\pm 1}$ for $i\in[N]$ is free abelian of rank $N$.  Denote its group algebra over $\ttk$ by 
$\H$. 
\itemb We have a  triangular decomposition,
$$U_q(\fg) = U_q(\fn^-) \ot \H \ot U_q(\fn^+).$$
 \itemc The PBW theorem for $U_q(\fn^-)$: the elements 
$f_{m_1, n_1}  \ldots  f_{m_s, n_s}$ with $m_i, n_i \in [N+1]$, $m_i < n_i$ and 
$(m_1, n_1)\le   \ldots \le  (m_s, n_s)$ form a basis for $U_q(\fn^-)$.
\ei
\et \bpf See \cite{Y}.\epf \noi 
If $\{\epsilon_i\mid i\in[N+1]\}$ is the standard basis for the dual of the Cartan subalgebra of $\fgl(N+1)$  consisting of diagonal matrices,  
then the simple roots of $\fg$ have the form
\be \label{dig}\ga_i = \gep_{i} -\gep_{i+1},\ee
 for $i\in [N]$. Let $Q = \sum_{i=1}^N \Z \ga_i$ and 
$Q^+ = \sum_{i=1}^N \N \ga_i$ respectively denote the root lattice and the positive root lattice of $\fg$. In addition set $\gs_i = \gep_{1} -\gep_{i+1}.$
There is an isomorphism $Q\lra \A$ sending $\ga_i$ to $k_i$. 
If $\ga\in Q,$ let $k_\ga \in \A$ be the image of $\ga$ under this map.
Set $U_q(\fb^\pm) = U_q(\fn^\pm) \ot \H$. This definition does not make  it clear that~$U_q(\fb^\pm)$ is a subalgebra of $ U_q(\fg)$.  However identifying $x\ot y$ with $xy$, we see that~$U_q(\fb^\pm)$ identifies with the product (not a direct product) of algebras $U_q(\fn^\pm) \H$.  Since the 
$k_\ga$ with $\ga \in \A$ form a $\ttk$ basis for $\H$,
every element of $ U_q(\fn^\pm)\H$ can be written in a unique way as a finite sum  $\sum_{\ga\in \A} a_\ga k_\ga$, where 
$a_\ga\in U_q(\fn^\pm)$. Equation~\eqref{eq1.1} implies that 
$U_q(\fn^\pm) \H = \H U_q(\fn^\pm),$ and from this it follows that $U_q(\fb^\pm)$ is a subalgebra. 
\\ \\
  The ideal of $U_q(\fn^+)$ generated by the $e_i$ is denoted $I{(\fn^+)}$.  
Since $ \H$ is free over $\ttk$, $\_\_\ot \H$ is an exact functor from left 
$U_q(\fn^+)$-modules  to left 
$U_q(\fb^+)$-modules. Applying this functor to the exact sequence 
$$ 0 \lra
I{(\fn^+)} \lra U_q(\fn^+)\lra \ttk \lra 0,$$
we obtain
$$ 0 \lra
I{(\fn^+)}\ot\H \lra U_q(\fb^+)\lra \H \lra 0.$$
Because of  \eqref{eq1.1} we see that 
$I{(\fn^+)}\ot \H = \H I{(\fn^+)}$ 
is a two-sided ideal of 
$U_q(\fb^+)$, so 
$U_q(\fb^+)/(I{(\fn^+)}\ot\H)\cong \H$ as $\ttk$-algebras. This means that any $\H$-module can be regarded as a $U_q(\fb^+)$-module with $\H I{(\fn^+)}$ contained in its annihilator. 
\\ \\
We mainly use Yamane's definition because it is easy to work  with his PBW basis.  However for representation theory we need to  compare it to the definition in Jantzen's book \cite{J} Chapter 4. 
 In Type A (more generally in the simply laced case), there is a bilinear form $(\;,\;)$ on $Q$ defined by $(\ga_i,\ga_j)=a_{i,j}$ where $A=(a_{i,j})$ is  the Cartan matrix.  
Since  any root $\alpha $ is conjugate to a simple root under the Weyl group $W$, and the form  $(\;,\;)$ is known to be $W$-invariant, $\ga$ satisfies $(\alpha,\alpha )=2$, so that $\alpha=\alpha^{\vee}$ and $q_{\alpha}=q$ in \cite{J} 4.2 (1).
 Then in \cite{J}, $U_v(\frak{g})$ is the $\ttf(v)$-algebra with generators $E_i, F_i,  K_i^{\pm 1}$ for $i\in[N]$, subject to the relations \cite{J}, 4.3 (R2)-(R4).
\be \label{R3} K_iK_j= K_jK_i,\ee 
\be \label{R23} K_iE_jK_i^{-1}=v^{a_{i,j}}E_j,  \;\;K_iF_jK_i^{-1}=v^{-a_{i,j}}F_j, \ee
\begin{equation}\label{R4}
E_iF_j-F_jE_i=\delta_{ij}(K_i-K_i^{-1})/(v-v^{-1}).
\end{equation}
and the Serre relations \cite{J}, 4.3 (R5) and (R6). For example (R6) says that for $i\neq j$, we have
\be  \label{R6}  \sum_{i=0}^{1 -a_{ij}}(-1)^i\left[ \begin{array}{c}
                1 -a_{ij}\\i
                 \end{array}\right]_v  F_i^{1 -a_{ij}-i} F_j F_i^{i},\ee
and (R5) is an analogous relation where $F_k$ is replaced by $E_k$. 
Serre relations are so called because they occur in Serre's theorem giving generators and relations for a simple Lie algebra, \cite{H1} Proposition 18.1.

Now comparing with \eqref{eq1.1} and \eqref{owl}, and checking the compatibility of the Serre relations, 
 shows that when $v=q^2$, there is an algebra map
$U_v(\frak{g}) \lra U_q(\frak{g})$ given by 
\be \label{arm} E_i \lra e_i,\,\, F_i \lra f_i,\,\,K_i \lra k_i^2.\ee
Many  of our computations take place in 
$U_v(\frak{g})$, so we identify 
$U_v(\frak{g})$ with its image under the above map.
\\ \\
The algebra $U_q(\fg) $ is $Q$-graded by setting $\deg e_i =- \deg f_i = \ga_i,$ and $\deg k_i=0$ for all $i$.  For $j>i+1$ we have
\be \label{cat0}  \deg f_{i,j}= -\sum_{k=i}^{j-1} \ga_k.
\ee   If $x\in U_q(\fg) $ and $\deg x = \gb,$ then 
\be\label{cat} k_\ga x k_\ga^{-1} = q^{(\ga,\gb)}x,\ee  compare \cite{J} 4.7 (1).
 In particular this gives an easy way to handle commutation relations between the~$k_\ga$ and the $f_{i,j}$. If $x\in U_q(\fg) $ has $\deg x = \gb,$ we write $x \in U_q(\fg)^\gb $. Set 
$K_\ga=k_\ga^{2}$ and $v=q^2$. Then \eqref{cat} becomes 
\be\label{cat1}K_\ga x K_\ga^{-1} = v^{(\ga,\gb)}x.\ee

\subsection{Representation Theory.} From now on we assume  $v=q^2$.
Introduce the {\it weight lattice } $P$ by 
$$P= \{\gl\in\fh^*| (\gl,Q) \subseteq \Z\}.$$  Since the form $(\;,\;)$ is integer valued on $Q$ we have $Q \subset P$.  In fact $|P/Q| = N+1$, the determinant of the Cartan matrix. 
If $M$ is a 
$U_v(\frak{g})$-module (or even an $\H$-module)  and $\gl \in P$, we say that $m\in M$ has {\it weight} $\gl$ if
 \be\label{cat3}K_\mu m = v^{(\gl,\mu)}m,\ee 
for all  $\mu \in Q$, see \cite{J} 5.1 (1). We consider only weight modules of type 1, see \cite{J} 5.2. Because of \eqref{arm}, this means that if 
 $M$ is an
$U_q(\frak{g})$-module  and $\gl \in P$, we say that $m= m_\gl \in M$ {\it has weight} $\gl$ if $k_\mu m = q^{(\gl,\mu)}m$ for all  $\mu \in Q$.
Let $\ttk v_\gl$ be the one dimensional $\H$-module with weight $\gl\in P$.  By the above remarks, we can regard $\ttk v_\gl$ as an $
U_q(\fb^+) 
$-module annihilated by $I{(\fn^+)}\ot \H$. Then define the {\it Verma module} $M(\gl)$ with highest weight $\gl$, for $U_q(\frak{g})$ 
by $M(\gl)=U_q(\frak{g}) \ot_{U_q(\fb^+) } \ttk v_\gl.$  The set of weights of $M(\gl)$ is $\gl - Q^+.$
\\ \\
If $r \in \N$ set $[r]_v = (v^{r}-v^{-r})/(v-v^{-1}).$ 
We use  some notation,
 $$[K_\mu;a]=(K_\mu v^a-K_\mu^{-1}v^{-a})/(v-v^{-1}),$$
 from \cite{J} 1.3 (1). 
Note that \eqref{cat3} 
 implies 
\be\label{cat7} [K_{\mu}:a]m_\gl = [a+ (\gl,\mu)]_v m_\gl.
\ee
\subsection{ Shapovalov elements.} \label{se}

Let $\fg$ be a simple Lie algebra with set of simple and positive roots $\Pi$ and $\Delta^+$ respectively. We denote the half sum of  positive roots by $\gr$. The
definition and significance of Shapovalov elements is given in \cite{M2} Section 1. We adapt the definitions  to the quantum case. For ease of exposition, and because we will only work with Type A in the sequel, we assume that $\fg$ is simply laced, and as before all roots $\alpha $ satisfy $(\alpha,\alpha )=2$, so  $\alpha=\alpha^{\vee}$.
\\ \\ 
Fix a positive root $\eta$  and a  positive   integer $m$. 
Let $P(m\eta)$ be the set of partitions of $m\eta$, see
\cite{M} Chapter  8  for notation. Then let
$\pi^0 \in P(m\eta)$ be the unique partition of $m\eta$ such
that $\pi^0(\ga) = 0$ if $\ga \in \Delta^+ \backslash \Pi.$
We say that $\gth = \theta_{\eta,m}\in U_q({\mathfrak b}^{-} )^{- m\eta}$ is a
{\it Shapovalov element for the pair} $(\eta,m)$ if it has the form 
\be \label{rat}
\theta = \sum_{\pi \in P(m\eta)} e_{-\pi}
H_{\pi},\ee where $H_{\pi} \in  \H $, $H_{\pi^0} = 1,$  and
\be \label{boo} e_{\ga} \theta \in 
U_q({\mathfrak g})(q^{(\gr, {\eta})}k_{\eta}-q^m)
+U_q({\mathfrak g})I{(\fn^+)} , \; \rm{ for \; all }\;\ga \in \Pi.\ee
For a semisimple Lie algebra, the existence of such elements was shown by Shapovalov, \cite{Sh} Lemma 1. Shapovalov elements for the pair $(\eta,m)$ are only unique up to the addition of an element from the left ideal in $U_q(\fg) $ appearing on the right side of \eqref{boo}, this is the quantum analog of  \cite{M1} Theorem 2.1. For any element $\gth$ as in \eqref{rat}  and $\gl\in P$, define the evaluation of $\gth$ at  $\gl\in P$ to be 
$\theta(\gl) = \sum_{\pi \in P(m\eta)} e_{-\pi}
H_{\pi}(\gl).$  Note that $\theta v_\gl = \theta(\gl) v_\gl.$ 
\\ \\
Consider the hyperplane in ${\mathfrak h}^*$ given by
 \be \label{vat}{\mathcal H}_{\eta, m} = \{ \lambda \in  {\mathfrak h}^*\mid (\lambda + \rho, \eta) =m
 \}.\ee
From \eqref{boo}, the Shapovalov element
$\theta_{\eta,m}$
has the important property that if $\gl
\in P\cap {\mathcal H}_{\eta, m}$ then
$\theta_{\eta, m}v_{\gl}$  is a highest weight vector of weight $\lambda -m\eta$ in
$M(\gl)$.
The normalization condition $H_{\pi^0}=1$
guarantees that $\theta_{\eta, m}v_\gl$ is never zero. We give an explicit expression for a Shapovalov elements for $U_q(\fsl(N+1))$. There are two important reductions we can make. First we need only consider the highest root  $\eta = \gep_1 -\gep_{N+1}$ of  $\fsl(N+1)$,  because any positive root is the highest root of some special linear subalgebra.
The second is somewhat deeper;
by a careful consideration of powers of $\gth_{\eta,1}$
in Section \ref{pse} it is only  necessary to consider the case $m=1$. So to simplify notation we set $\mathcal{H}_{\eta} = \mathcal{H}_{\eta,1},$ and
denote a  Shapovalov element for the pair $(\eta, 1)$ by $\gth_{\eta}$.

\section{Commutation Relations.}\label{sec2}
\noi 
We need several commutation relations 
for the computation of the Shapovalov elements. These may also be found in \cite{Y}. For elements $a, b$ in an associative algebra, set $[a, b]=ab-ba$. We use the elementary identity $[a, bc]=b[a, c]+[a, b]c$, many times without comment. 
For comparison with the relations in \cite{Y}, set $e_{i,i+1}= e_i$. 
In particular, we consider the question of when $[e_\ell, f_{i,j}] \neq 0.$ If $j = i+1$, this occurs if and only if $i = \ell$ by \eqref{owl}.  Now suppose that $j > i+1$.
\bl\label{sow} If $j > i+1$ and $[e_\ell, f_{i,j}] \neq 0$, then $\ell = i$ or 
$\ell +1 = j$.\el
\bpf  Use induction on $j-i-1$, and Equations \eqref{owl} and \eqref{dog}. Assume  $\ell \neq  i$ and $\ell \neq  j-1$, and show $[e_\ell, f_{i,j}] = 0$.
By definition $f_{i,j}=qf_{i,j-1}f_{j-1}-q^{-1}f_{j-1}f_{i,j-1}.$ 
Since  $\ell\neq j-1$, $[e_\ell, f_{j-1}]=0$. If $j-1 =i+1$, then 
$[e_\ell, f_{i,j-1}] =[e_\ell, f_{i}] $,  and this is zero since 
 $\ell \neq  i$. Otherwise 
$j-1 >i+1$ and 
$[e_\ell, f_{i,j-1}] =0$ by induction.
 \epf
\noi In the exceptional cases $\ell = i$  (resp. $\ell+1 = j$) 
the situation is covered by the next two Lemmas, see \eqref{eq2.2}, (resp.  \eqref{eq2.3}) with some small changes of notation which will be useful later. In both proofs we make use  of \eqref{cat}. 
The relations are  special cases of the first (resp. sixth) equation in Section 3 (3)  of \cite{Y}.  Since these relations are crucial to our work, but  not proved in \cite{Y} we give the proofs.
\begin{lemma}
When $i\leq b-2$,
\begin{equation}\label{eq2.2}
[e_{i},f_{i,b}]=qf_{i+1,b}k_i^2.
\end{equation}
\end{lemma}

\begin{proof} We use induction on $b$.
\newline
{Case 1}: If $b=i+2$, then $[e_{i},f_{i+1,i+2}] =0$ by \eqref{owl}. Therefore using \eqref{owl} again,
\[
\begin{aligned}
&[e_{i},f_{i,i+2}]=[e_{i},qf_{i,i+1}f_{i+1,i+2}-q^{-1}f_{i+1,i+2}f_{i,i+1}]\\
&=(q^2-q^{-2})^{-1}[(k_i^2-k_i^{-2})\,q\,f_{i+1,i+2}-q^{-1}f_{i+1,i+2}(k_i^2-k_i^{-2})]\\
&=(q^2-q^{-2})^{-1}[f_{i+1,i+2}\,q\,(q^2k_i^2-q^{-2}k_i^{-2})-q^{-1}f_{i+1,i+2}(k_i^2-k_i^{-2})]\\
&=(q^2-q^{-2})^{-1}[q^3f_{i+1,i+2}k_i^2-q^{-1}f_{i+1,i+2}k_i^2]\\
&=qf_{i+1,i+2}k_i^2.
\end{aligned}
\]
{Case 2}: If $b>i+2$, then $[e_{i}f_{b-1,b}]=0$  by \eqref{owl}, and by induction we know the result for $[e_{i},f_{i,b-1}]$. Thus,
\[
\begin{aligned}
&{[e_{i},f_{i,b}]}=[e_{i},qf_{i,b-1}f_{b-1,b}-q^{-1}f_{b-1,b}f_{i,b-1}]\\
&=q^2\,f_{i+1,b-1}k_i^2f_{b-1,b}-q^{-1}f_{b-1,b}\,q\,f_{i+1,b}\,k_i^2\\
&=q\left(qf_{i+1,b-1}f_{b-1,b}-q^{-1}f_{b-1,b}f_{i+1,b}\right)k_i^2\\
&=qf_{i+1,b}k_i^2.
\end{aligned}
\]
\end{proof}

\begin{lemma}
When $a<i$,
\begin{equation}\label{eq2.3}
[e_{i},f_{a,i+1}]=-q^{-1}f_{a,i}k_i^{-2}.
\end{equation}
\end{lemma}

\begin{proof}
By definition,
\[
f_{a,i+1}=qf_{a,i}f_{i,i+1}-q^{-1}f_{i,i+1}f_{a,i},
\]
and clearly
\[
[e_{i},f_{a,i}]=0,
\]
due to Lemma \ref{sow} when $a<i-1$ and to \eqref{owl} when $a=i-1$. So by \eqref{owl}, 
\[
\begin{aligned}
&{[e_{i},f_{a,i+1}]}=(q^2-q^{-2})^{-1}[q\,f_{a,i}(k_i^2-k_i^{-2})-q^{-1}(k_i^2-k_i^{-2})f_{a,i}]\\
&=(q^2-q^{-2})^{-1}[q\,f_{a,i}(k_i^2-k_i^{-2})-q^{-1}f_{a,i}(q^2k_i^2-q^{-2}k_i^{-2})]\\
&=(q^2-q^{-2})^{-1}f_{a,i}(-q+q^{-3})k_i^{-2}\\
&=-q^{-1}f_{a,i}k_i^{-2}.
\end{aligned}
\]
\end{proof}\noi
Consider a strictly increasing sequence of integers
\be \label{dod} I = \{j_0 , j_1, \ldots ,j_s, j_{s+1} \},\ee
with $1\le j_0$ and $j_{s+1} \le N+1$.
If $I$ is a singleton set, then by definition set $f_I =1$. Otherwise define  $f_{I} \in U_q(\fn^-)$ by
 \by \label{ddd}f_I = f_{j_{0},j_1}f_{j_1,j_2} f_{j_2,j_3} \ldots f_{j_s,j_{s+1}}. \ey
The way the factors are ordered in \eqref{ddd} is consistent with the lexicographic order used in Yamane's PBW Theorem, Theorem \ref{YT}.   This is because   $j_{k-1} < j_{k} < j_{k+1}$
implies that $(j_{k-1},j_{k}) < (j_{k}, j_{k+1})$ in the lexicographic order.
Let
\be \label{120} \I = \{I\subseteq [N+1]\mid 1, N+1\in I\}, \ee  
and for $ I\in \I$, define
\be \label{130} r(I) = \{s-1\mid s \in \bar I\} \ee
where $\bar I$ be the complement of $I$ in $\I$.
Note that if $ I\in \I$, then by \eqref{ddd} and the definition of $\I,$ $f_I \in U_q(\fn^-)^{-\eta}$.  Clearly elements of $\I$ correspond bijectively  to partitions of $\eta$, (with a suitable ordering on positive roots) and taking $J=\I$ we have $f_J = e_{-\pi^0}.$ Thus the $f_I$ with $I\in\I,$ form  a $\ttk$-basis for $
U_q(\fn^-)^{-\eta}$.  
\\ \\
Consider the set 
\be \label{140} \cS_i=\{I\subseteq \I\mid i,i+1 \in I   \}.\ee
\noi For $I\in \cS_i$, define 
\be\label{qrt} I^+ = I \backslash \{i\}, \mbox{ and } I^- = I \backslash \{i+1\}.\ee
If $i=N$, then $I^- \notin \I$ and if $i=1,$ then $I^+ \notin \I$. We set
$f_{I^-}=0$ or $f_{I^+}=0$ in these cases. 
It is often useful to write $I\in \cS_i$ in the form
\be \label{az1} I =\{1, \ldots ,a,i,i+1,b, \ldots , N+1 \}.\ee 
Then 
\be \label{az2} I^+ =\{1, \ldots ,a,i+1,b, \ldots , N+1 \}\mbox{ and }
I^- =\{1, \ldots ,a,i,b, \ldots , N+1 \}.\ee 
The following result summarizes the  commutation relations that we need later in the paper.
\bl \label{bat} If $ J\in \I$, then $e_i$ fails to commute with at most one factor of~$f_J$.  Furthermore, for a fixed $I\in \cS_i$, we have:
\bi \itema $f_i$ is a factor of $f_{I},$ and  $[e_i, f_i] $ is given by \eqref{owl}.
\itemb 
$f_{i,b}$ is a factor of $f_{I^-},$ and   $[e_{i},f_{i,b}]$ is given by \eqref{eq2.2}.
\itemc 
$f_{a,i+1}$ is a factor of $f_{I^+}, $and $[e_{i},f_{a,i+1}]$ is given by \eqref{eq2.3}.
\itemd $e_i$ commutes with all other factors of 
$f_{I}$ and  $f_{I^\pm}$ 
in parts $(a)$-$(c)$.
\iteme If $J$ does not have the form $I, I^+$ or $I^-$ for some  $I\in \cS_i$, then $e_i$ commutes with all factors of $f_J$.
\ei
\el  

\bpf The claims about the factors in (a)-(c) hold by \eqref{az1} and \eqref{az2}. The rest follows from Lemma \ref{sow}.\epf
\section{Shapovalov elements for the quantized enveloping algebra $U_q(\fsl(N+1))$} \label{sec3} \noi
 For $i\in[N]$ and $J\in \mathbb{I}$, define
\be \label{1q}  h_i =-q^{-1}
v^{-((\gr,\gs_i)-1)} K^{-1}_{{\gs_{i}}}[K_{\gs_i}:i]_v \in \H,
\quad H_J =\prod_{i \in r(J)}h_{_i}.\ee
For the rest of the paper, $v_{\lambda}$ is a highest weight vector of weight $\lambda $ in the Verma module $M(\lambda )$.
\bl We have
\be \label{0q} h_iv_\gl=-q^{-1}v^{-((\gl +\gr,\gs_i)-1)}[(\gl+\gr,\gs_i)]_v v_\gl\ee for all $\gl \in Q.$
\el
\bpf
Since $\gs_i$ is the sum of $i$ simple roots, and for a simple root $\ga,$  we have
$(\gr, \ga) =(\gr, \ga ^\vee)= 1,$ it follows that  $(\gr, \gs_i)=i$.  Thus 
 $(\gl+\gr,\gs_i)= (\gl,\gs_i) +i$, so by \eqref{cat7},  
\be \label{bo1}[K_{\gs_i}:i] v_\gl = [(\gl+\gr,\gs_i)]_v v_\gl.\ee
Also
\be \label{bo2} v^{-((\gr,\gs_i)-1)} K^{-1}_{{\gs_{i}}} v_\gl =v^{-((\gl +\gr,\gs_i)-1)}  v_\gl.
\ee
The result follows by combining \eqref{bo1}
and \eqref{bo2}.
\epf \noi 
Our first main result provides an  explicit expression for a Shapovalov element~$\gth_\eta$.

\bt\label{bb} Let $\eta =\gep_1 -\gep_{N+1}$, and suppose $v_{\lambda }$ is a highest weight vector in the $\fg$-Verma module $M(\gl)$ of weight $\gl$, and set
\be \label{cat22} \Gt_\eta= \sum_{J\subseteq \I} f_J H_J.\ee
Then \bi \itema $e_{\ga_k}\Gt_\eta v_{\lambda }=0$ for $k\in [N-1]$.
\itemb If $(\gl+\gr,\eta)=1,$  then $\Gt_\eta v_{\lambda }$ is a $\fg$-highest weight vector, and so
$\Gt_\eta$ is a Shapovalov element for the pair $(\eta,1)$.\ei
\et  
\noi By the remarks preceding \eqref{140} we see that the coefficient of 
$ e_{-\pi^0}$ in  \eqref{cat22}  is $H_{\pi^0} =1.$  Thus \eqref{rat} holds.
We show that $\sum_{J\subseteq \I} f_J H_J$ satisfies Equation \eqref{boo}. 
Given a simple root $\ga_i=\gep_{i} - \gep_{i+1}$, to avoid double subscripts we sometimes set $\ga= \ga_i$.
Thus  $e_\ga= e_{i}$. 
The goal is to show that if $\gl\in\cH_\eta$, then  

\be \label{cat44}
e_i \sum_{J\subseteq \I} f_J H_J v_{\lambda }=0
\ee for all $i\in [N]$. 
This will be shown directly in this Section.  In Subsection \ref{dxd} we give a different approach to Shapovalov elements and it is convenient to make some remarks at this point that apply to the second approach. 
We can rewrite Equations~\eqref{cat44}  as a sequence of equalities in $\H$, depending on $\gl$, one for each partition of $\eta-\ga_i$.  It is enough to show that the equalities hold for all $\gl$ in a Zariski dense subset $\gL$ of  $\cH_\eta$, with  $\gL  \subset P$. The subset $\gL$ will be defined in \eqref{1tar}. 
The second approach can be used to construct Shapovalov elements $\theta_{\eta, m} $ for $m>1$.   
The two approaches are compared in  Subsection \ref{2cal}.  
\\ \\
We conclude the proof of  
Theorem  \ref{bb} after a series of Lemmas.
\\ \\
\noi
Now fix $i$ and recall  the set  $\cS_i$ from \eqref{140}. If $I \in \cS_i$ then  $e_{-\ga} = f_i$ is a factor of $f_I.$
With $I$ as in \eqref{az1}, 
we set  
\be \label{ay} I_1 = \{1, \ldots ,a,i \}\mbox{ and } I_2 = \{i+1, b, \ldots , N+1 \}. \ee
Then   we have


\bl \label{cow} 
If $v_{\lambda }$ is a highest weight vector of weight $\lambda $, then 
\by 
 e_i f_{I}v_{\lambda } &=& (q^2-q^{-2})^{-1}f_{I_1} (k_i^2-k_i^{-2})f_{I_2}v_{\lambda }\nn\\& =&
\left\{ \begin{array}
  {cc}(q^2-q^{-2})^{-1}f_{I_1}(k_i^2-k_i^{-2}) v_{\lambda }&\mbox{if} \;\; I_2  \mbox{ is a singleton }\label{2vam}
\\ (q^2-q^{-2})^{-1}f_{I_1}f_{I_2}(q^2k_i^2-q^{-2}k_i^{-2}) v_{\lambda } &\mbox{ otherwise}.
\end{array} \right. \ey \el

\bpf Note that $f_{I} = f_{I_1}f_if_{I_2}.$ We have $[e_i, f_{I_1}]=0$ by 
Lemma  \ref{bat}. The first equality follows from this and   \eqref{owl}. If $|I_2|>1$, the first factor in $f_{I_2}$ has the form $f_{i+1,b}$ and by \eqref{cat}
$$(k_i^2-k_i^{-2})f_{i+1,b}=   f_{i+1,b}(q^2k_i^2-q^{-2}k_i^{-2}).$$  It is easy to see that $k_i$ commutes with the other 
factors of  $f_{I_2}.$ 
\epf

\noi
Next with $I$ as in \eqref{az1}, set \be \label{pig}I_1^+ = \{1, \ldots ,a\}\mbox{ and } I_2^- = \{b, \ldots , N+1 \}, \ee
and note the factorizations
\[  f_{I^+} = f_{I_1^+}f_{a,i+1}f_{I_2}\mbox{ and } f_{I^-}= f_{I_1}f_{i,b}f_{I_2^-}.
\] 
If $v_{\lambda }$ is a highest weight vector of weight $\lambda $, we also need to compute 
$[e_i, f_{I^\pm}]v_{\lambda }$. Since 
$e_iv_{\lambda } =0$, this is equivalent to finding $[e_i, f_{I^\pm}]$. Now using Lemma \ref{bat} and  \eqref{cat}, we have
\by \label{asp} [e_i, f_{I^+}]&=& 
[e_i, f_{I_1^+}f_{a,i+1}f_{I_2}]  = f_{I_1^+}[e_i, f_{a,i+1}]f_{I_2}\\
&=&-q^{-1} f_{I_1^+} f_{a,i}k_i^{-2}f_{I_2}\nn\\
&=& 
\left\{ \begin{array}
  {cc}-q^{-1} f_{I_1} f_{I_2} k_i^{-2} &\mbox{if} \;\; I_2  \mbox{ is a singleton }
\\ -q^{-3} f_{I_1} f_{I_2} k_i^{-2}  &\mbox{ otherwise}.
\end{array} \right. \nn
\ey Note that $i<N$ implies $|I_2|\neq 1$ in 
\eqref{2vam} and 
\eqref{asp}.
Also by 
\eqref{eq2.2} and \eqref{cat},
\by \label{ox} [e_i, f_{I^-}] &=&  [e_i, f_{I_1}f_{i,b}f_{I_2^-}] =f_{I_1}[e_i,f_{i,b}]f_{I_2^-} \nn\\
&=&
qf_{I_1}f_{i+1,b}k_i^{2}f_{I_2^-}  \nn\\
&=&qf_{I_1}f_{I_2} k_i^{2}.
\ey 
We record an easy lemma.
\bl \label{hog} We have \bi 
\itema $r(I^+)=r(I)\cup \{i-1\} \mbox{ and } r(I^-)=r(I)\cup \{i\}$.
\itemb $H_{I^+} = h_{i-1}H_I$ and $H_{I^-} = h_{i}H_I.$
\itemc If $J\subseteq \I$ and $e_i f_{J}v_{\lambda }$ is a non-zero element of $ f_{I_1} f_{I_2}\H v_{\lambda }$, then
 $$J=I, I^+\mbox{ or } I^-.$$ If $I^+$ or  $I^-$ is not a subset of $\I$, they should be left out in the above equation.
\ei \el
\bpf (a) follows from the definitions \eqref{130} and \eqref{qrt}. Then (b) follows from \eqref{1q}.  For (c) note that the pair $(I_1, I_2)$ determines $I$ (and hence also $I^\pm$), since $I$ is the  disjoint union $I = I_1 \cup I_2$. By definition $I$ and $i$ determine $I_1$ and $I_2$.
\epf
\noi 
We begin with the case $i=1$.  Recall that $f_{I^+}=0$ in this case. 
\bl \label{Lem1}
$e_1(f_IH_I+f_{I^-}H_{I^-})v_{\gl }=0$.
\el
\bpf If $i=1$,  apart from some elementary simplifications, we use in order, 
 Lemma~\ref{hog}~(b),  \eqref{ox}, \eqref{0q}, and  Lemma~\ref{cow},
\by \label{as}e_1f_{I^-}H_{I^-}v_{\gl }&=& e_1f_{I^-}h_1H_{I}v_{\gl }
\nn\\&=& -f_{I_1}f_{I_2}v^{-((\gl +\gr,\ga_1)-1)}[(\gl+\gr,\ga_1)]_v  K_{\ga_1}H_{I}v_\gl
\nn\\&=& -f_{I_1}f_{I_2}v^{-(\gl ,\ga_1)} 
v^{(\gl,\ga_1)}[(\gl+\gr,\ga_1)]_v  H_{I}v_\gl
\nn\\&=&
-[(\gl+\gr,\ga_1)]_vf_{I_1}f_{I_2}H_Iv_{\gl }=-e_1f_IH_Iv_{\gl }.\ey 
\epf  \noi
\bl \label{hot}
\bi \itema If $1 <i \le N$, then 
\by \label{equ4}
 e_if_{I}H_{I}v_{\gl }&=& 
\left\{ \begin{array}
  {cc}f_{I_1}f_{I_2}[(\gl,\ga_i)]_vH_Iv_{\gl},&\mbox{if} \;\; 
i=N,
\\ f_{I_1}f_{I_2}[1+(\gl,\ga_i)]_vH_Iv_{\gl},&\mbox{ otherwise}.
\end{array} \right. \ey 
\itemb If $1 <i \le N$, then 
 \by \label{equ3}\;\;
e_if_{I^+}H_{I^+}v_{\gl }&=& 
\left\{ \begin{array}
  {cc}
q^2f_{I_1}f_{I_2}v^{-(\gl+\gr,\gs_i)}[(\gl+\gr,\gs_{i-1})]_vH_{I}v_{\gl },
&\mbox{if} \; i=N,
\\ 
f_{I_1}f_{I_2}v^{-(\gl+\gr,\gs_i)}[(\gl+\gr,\gs_{i-1})]_vH_{I}v_{\gl },
&\mbox{otherwise.}
\end{array} \right. \ey 
\itemc
If $1 <i<N$, then
\be \label{equ2}\begin{aligned} e_if_{I^-}H_{I^-}v_{\gl }&=-v^{-(\gl+\gr,\gs_{i-1})}[(\gl+\gr,\gs_i)]_vf_{I_1}f_{I_2}H_Iv_{\gl }.
\end{aligned}
\ee \noi
\ei
\el
\bpf (a) This follows from Lemma~\ref{cow} and \eqref{cat7}.\\ \\
(b) Suppose first that  $1 <i < N$.   Then by Lemma \ref{hog} (b), (\ref{asp}) and (\ref{0q}),
\be 
\begin{aligned}
e_if_{I^+}H_{I^+}v_{\gl }&=e_if_{I^+}h_{i-1}H_{I}v_{\gl }=f_{I_1}f_{I_2}q^{-3}h_{i-1}H_{I}k_i^{-2}v_{\gl }\\
&=q^{-4}f_{I_1}f_{I_2}H_{I}k_i^{-2}v^{-((\gl+\gr,\gs_{i-1})-1)}[(\gl+\gr,\gs_{i-1})]_v v_{\gl }\\
&= f_{I_1}f_{I_2}H_{I}v^{-(\gl+\gr,\ga_i)
}v^{-(\gl+\gr,\gs_{i-1})}[(\gl+\gr,\gs_{i-1})]_v v_{\gl }
\\
&=f_{I_1}f_{I_2}v^{-(\gl+\gr,\gs_i)}[(\gl+\gr,\gs_{i-1})]_vH_{I}v_{\gl }.\nn
\end{aligned}
\ee 
The proof is similar if $i=N$ except that since a different case holds in (\ref{asp}),
the first line of the proof becomes 
$$
e_if_{I^+}H_{I^+}v_{\gl }=f_{I_1}f_{I_2}q^{-1}h_{i-1}H_{I}k_i^{-2}v_{\gl },$$ leading to the result stated.
\\ \\
(c)
By Lemma \ref{hog} (b), (\ref{ox}) and (\ref{0q}),
\be \begin{aligned} e_if_{I^-}H_{I^-}v_{\gl }&=e_if_{I^-}h_iH_{I}v_{\gl }
\\
&= -f_{I_1}f_{I_2}v^{-((\gl +\gr,\gs_i)-1)}[(\gl+\gr,\gs_i)]_v  K_{\ga_i}H_{I}v_\gl
\\
&= -f_{I_1}f_{I_2}v^{-(\gl +\gr,\gs_i)}[(\gl+\gr,\gs_i)]_v  v^{(\gl +\gr,\ga_i)}H_{I}v_\gl
\\
&=-v^{-((\gl+\gr,\gs_i)-(\gl+\gr,\ga_i))}[(\gl+\gr,\gs_i)]_vf_{I_1}f_{I_2}H_Iv_{\gl }\\
&=-v^{-(\gl+\gr,\gs_{i-1})}[(\gl+\gr,\gs_i)]_vf_{I_1}f_{I_2}H_Iv_{\gl }.\nn
\end{aligned}
\ee 
\epf \noi
Until further notice we assume that $1 <i<N$. 
\bl \label{Lem2}
$$[(\gl+\gr,\ga_i)]_v+v^{-(\gl+\gr,\gs_{i})}[(\gl+\gr,\gs_{i-1})]_v-v^{-(\gl+\gr,\gs_{i-1})}[(\gl+\gr,\gs_i)]_v=0.$$
\el
\bpf  Using $\gs_{i-1}-\gs_{i}=
-\ga_{i}$ to cancel the second and third terms in the first line below, we have
$$\begin{aligned}
&\left(v^{(\gl+\gr,\ga_{i})}-v^{-(\gl+\gr,\ga_{i})}\right)+v^{-(\gl+\gr,\gs_{i})}\left(v^{(\gl+\gr,\gs_{i-1})}-v^{-(\gl+\gr,\gs_{i-1})}\right)\\
&=v^{(\gl+\gr,\ga_{i})}-v^{-(\gl+\gr,\gs_{i})-(\gl+\gr,\gs_{i-1})}\\
&=v^{-(\gl+\gr,\gs_{i-1})}\left(v^{(\gl+\gr,\gs_{i})}-v^{-(\gl+\gr,\gs_{i})}\right).
\end{aligned}$$
This is equivalent to the statement of the lemma.
\epf

\bc \label{Cor3}
For $1<i<N$, we have
 $$e_i(f_IH_I+f_{I^+}H_{I^+}+f_{I^-}H_{I^-})v_{\gl}=0.$$
\ec

\bpf
We have using  \eqref{equ4}, \eqref{equ3} and \eqref{equ2},
\by \label{aa} &&e_i( f_{I}H_I + f_{I^+}H_{I^+} + f_{I^-}H_{I^-})v_{\gl} =\nn\\
&&  f_{I_1}f_{I_2}([1+(\gl,\ga_i)]_v+v^{-(\gl+\gr,\gs_i)}[(\gl+\gr,\gs_{i-1})]_v-v^{-(\gl+\gr,\gs_{i-1})}[(\gl+\gr,\gs_i)]_v
)H_Iv_{\gl}, \nn 
\ey
and this is zero by
Lemma~\ref{Lem2}. 
\epf

\begin{proof}[Proof of Theorem \ref{bb}]
(a) We first  prove the result for $1<i<N$.  
By Corollary~\ref{Cor3}, we have 
\begin{equation}\label{aa}
e_i( f_{I}H_I + f_{I^+}H_{I^+} + f_{I^-}H_{I^-}) v_{\gl} =0. 
\end{equation}
Now if $J$ does not contain $i$ or $i+1$, then
$e_i  f_J H_J v_{\lambda }=0$. Otherwise
there exists a set $I\in \cS_i$ 
 such that
$J = I, I^+$ or $I^- $. Thus we have using \eqref{aa},
\[ e_\ga \sum_{J\subseteq \I} f_J H_J v_{\lambda }= e_\ga \sum_{I \in \cS_k} (f_{I}H_I + f_{I^+}H_I^+ +f_{I^-}H_I^-)v_{\lambda } =0.\]
The proof for $i=1$ is entirely analogous. 
We use  Lemma \ref{Lem1} in place of  Corollary~\ref{Cor3}.
This proves the first statement in the Theorem. 

\medskip
(b) By   \eqref{equ4} and \eqref{equ3} we have  $e_N(f_IH_I+f_{I^+}H_{I^+})v_{\gl }=$ 
\[
\begin{aligned}
&=f_{I_1}f_{I_2}\{[(\lambda,\alpha_N)]_v+v^{1-(\lambda+\rho ,\sigma_N)}[(\lambda+\rho,\sigma_{N-1})]_v\}H_Iv_{\lambda }\\
&=f_{I_1}f_{I_2}\{v^{(\lambda,\alpha_N)}-v^{-(\lambda,\alpha_N)}+v^{1-(\lambda+\rho ,\sigma_N-\sigma_{N-1})}-v^{1-(\lambda+\rho ,\sigma_N+\sigma_{N-1})}\}
/(v-v^{-1})H_Iv_{\lambda }\\
&=f_{I_1}f_{I_2}\{v^{(\lambda,\alpha_N)}-v^{1-(\lambda+\rho ,2\eta -\alpha_N)}\}/(v-v^{-1})H_Iv_{\lambda },
\end{aligned}
\]
where the last equality follows since $\sigma_N-\sigma_{N-1}=\alpha_N$, $(\rho ,\alpha_N)=1$ and $\sigma_N=\eta $.  We conclude  the proof with the following Lemma.
\end{proof}

\begin{lemma}\label{L3.9}
$v^{(\lambda,\alpha_N)}=v^{1-(\lambda+\rho ,2\eta -\alpha_N)}$, if $\lambda \in H_{\eta }$.
\end{lemma}

\begin{proof} The difference of the two sides is
\be \label{R7}v^{(\lambda,\alpha_N)}-v^{1-(\lambda+\rho ,2\eta -\alpha_N)}=v^{(\lambda,\alpha_N)}(1-v^{1-2(\lambda +\rho ,\eta )+(\rho ,\alpha_N)})=0.\ee
The last equality follows from $(\lambda +\rho ,\eta )=1$, due to $\lambda \in H_{\eta }$, and $(\rho ,\alpha_N)=1$.
\end{proof}

\section{Powers of Shapovalov  elements.} \label{pse}
Powers of the Shapovalov  element $\theta_{\eta,1}$ for a simple Lie algebra $\fg$  of types A-D were considered in \cite{M1}. 
For $U_q(\fg)$ powers of $\theta_{\eta,1}$  were considered in  \cite{Mu1}, based on previous work by the same author on the R-matrix and inverse Shapovalov form, \cite{Mu}. In this Section we present an elementary approach to powers of 
$\theta_{\eta,1}$  in $U_q(\fsl(N+1))$ 
  by adapting the method of \cite{M1}. 
  It is based on the original method of Shapovalov, which is rather different from the method we have used so far in this paper.  Thus to tie up any loose ends we first indicate the connection between the two approaches.  It is convenient to do this using certain non-commutative determinants which we introduce in the next Subsection.

\subsection{Evaluation of Shapovalov element using noncommutative determinants.} 
We define a noncommutative determinant of the $n\times n$ matrix $B=(b_{ij})$,
working from left to right, by
\be \label{lr}{\stackrel{\longrightarrow }{{\rm det}}}(B) =  \sum_{w \in \mathcal{S}_n} \mbox{sign}(w) b_{w(1),1} \ldots b_{w(n),n}, \ee
Cofactor expansions of ${\stackrel{\longrightarrow }{{\rm det}}}(B) $ are valid as long as the overall order of the terms is unchanged. Next  let $m=N+1$, and for $i \in [m-2] $.  Define $h_i$, 
$H_J$ from \eqref{1q} and~\eqref{0q}, and  let $ c_i = h_i(\gl)$. Then set

\be\label{coc}
\cD^{N+1}(\gl)=  \left[ {\begin{array}{ccccc}
f_{1,2}&f_{1,3} & \hdots  & f_{1,N}  & f_{ 1,{N+1}} \\
 -c_{1} &f_{2, 3} & \hdots & f_{2, N}  & f_{2,{N+1} } \\
0& -c_{2} & \hdots& f_{3,N }&    f_{3,{N+1}} \\
 \vdots  & \vdots & \ddots & \vdots   & \vdots \\
0&0& \hdots & -c_{N -1} & f_{ N,{N+1}} 
\end{array}}
 \right].
\ee
\bt \label{T}
The Shapovalov element for $\eta$ satisfies
\[\gth_{\eta}(\gl) = {\stackrel{\longrightarrow}{{\rm det}}}\;\cD^{{N+1}}(\gl),\]
for all $\gl \in \cH_{\eta}$.
\et \noi
\bpf We show  
${\stackrel{\longrightarrow}{{\rm det}}}\;\cD^{{N+1}}(\gl)$
is the evaluation of 
the element $\Gt_{\eta}$ from Theorem~\ref{bb}
at $\gl\in \cH_{\eta}$. In other words, we show 
\be \label{da1}{\stackrel{\longrightarrow }{{\rm det}}}\;\cD^{{N+1}}(\gl)= \sum_{I\subseteq \I} f_I H_I(\gl).\ee Let $\I$ be  as  in \eqref{120}. 
We consider the  complete expansion of
${\stackrel{\longrightarrow}{{\rm det}}}\;\cD^{{N+1}}(\gl)$.
 Each term in the  expansion 
is obtained by choosing  a non-zero product of elements from each column, with each row occurring exactly once.   
The product of the 
 chosen elements lying  above the subdiagonal has the form  $f_I$ for some $I\in \I$. 
The proof of~\eqref{da1} is completed by the following Lemma.
\epf
\bl \label{ece} The product of subdiagonal terms accompanying $f_I$ is 
\be \label{dd} \prod_{i\in r(I)} c_{i} =\prod_{i\in r(I)} h_{i}(\gl) = H_I(\gl).\ee
 \el 

\bpf 
It is easy to adapt the proof from \cite{M2} Section 4. 
\epf  
\subsection{Some remarks on the adjoint action.} \label{dyd}
We make $U_v(\frak{g})$ 
 into a  Hopf algebra,  such that $U_v(\frak{b}^\pm)$ 
are Hopf subalgebras by defining the  coproduct $\Gd$ and antipode~$S$ on $U_v(\frak{b}^-)$ as in \cite{J} 9.13 (4)-(6) and (9).  In particular for the generators of~$U_v(\frak{b}^-)$, we have  
$$\Gd(F_\ga)  = F_\ga \ot1 + K_\ga \ot F_\ga, \quad \Gd(K_\ga)  = K_\ga\ot  K_\ga,$$ and $$S(F_\ga)  = -K_\ga^{-1} F_\ga , \quad S(K_\ga)  = K_\ga^{-1}.$$
\ff{For a different Hopf algebra structure on  $U_v(\frak{g})$, see  \cite{J} Lemma 4.8.
}
Hence by \eqref{arm} there is a Hopf algebra structure on
$ U_q(\frak{b}^-)$ given by
$$\Gd(f_\ga)  = f_\ga \ot1 + K_\ga \ot f_\ga, \quad \Gd(K_\ga)  = K_\ga\ot  K_\ga,$$ and $$S(f_\ga)  = -K_\ga^{-1} f_\ga ,\quad S(K_\ga)  = K_\ga^{-1}.$$
Now for any Hopf algebra $H$ and $e\in H$, write $\Gd(e) = \sum e_1\ot e_2$ in Sweedler notation.  Then  define the map $\ad_e:H\lra H$ by 
$\ad_e(x) = \sum e_1 x S(e_2)$. If 
$x\in H =  U_v(\frak{b}^-)$  and $F =F_\ga=f_\ga, K= K_\ga$ we have 
\be  \label{R8}\ad_F(x) = Fx - KxK^{-1}F = Fx - \gs(x)F\ee  where $\gs$ is the automorphism of given by
\be  \label{R9}\gs(x) = KxK^{-1}.\ee
Note that for $u \in U^\nu$ we have $\gs(u) = v^{(\ga,\nu)}u$ by \eqref{cat1}.
\\ \\
A short calculation (note that $1 -a_{\ga \gb} \le 2$)  shows that
if $\ga\neq\gb$, then 
\be  \label{R7} \ad_{F_\ga}^{ 1 -a_{\ga \gb}} (F_\gb)  =\sum_{i=0}^{1 -a_{\ga \gb}}(-1)^i\left[ \begin{array}{c}
                1 -a_{\ga \gb}\\i
                 \end{array}\right]_v  F_\ga^{1 -a_{\ga\gb}-\ga} F_\gb F_\ga^{i},\ee and this is zero by \eqref{R6}.
  Using \eqref{cat1} it is easy to check the following 
\bl\bi \itemo \itema
$\gd = \ad_{F_\ga}$  is  a $\gs$-derivation, that is 
$\gd (xy) = \gd(x)y +\gs(x)\gd(y)$.
\itemb $v^2( \gs\circ\gd )= \gd \circ \gs$. 
\ei\el \noi
Adapting the proof of  \cite{G} Lemma 6.2 which uses a different definition of Gaussian binomial coefficients, $\gd$ satisfies the Leibniz identity.

\be  \label{R24}  \gd^{n}(ab)= \sum_{i=0}^n v^{i(n-i)}\left[ \begin{array}{c}
                n \\
                i \end{array}\right]_v\gs^{i}( \gd^{n-i}(a))\gd^{i}(b).\ee
\subsubsection{The induction set-up.} \label{dcd}
Let $\ga_i$ be the simple roots of $\fg$ as in \eqref{dig}. Then consider the chain of subalgebras \be
\label{r29} \fk_1 \subset \fk_2 \ldots \subset \fk_N =\fg,\ee such that $\fk_i \cong \fsl(i+1)$. Let $\ga_1$ be the simple root of $\fk_1$,
and for $i>1$ let $\ga_i$ be the simple root of $\fk_i$ that is not a root  of  $\fk_{i-1}.$   
In the next Subsection, we show how to construct Shapovalov elements (with $m=1$ in \eqref{rat}) for the highest root of $\fk_i$ by induction on $i$.  We only need to show the last step of the induction, so to simplify we use the following notation for the rest of this Subsection.  Let  $\fk =  \fk_{N-1}.$  Then if  $\mathfrak{n}^+\oplus \mathfrak{h}\oplus \mathfrak{n}^-$ is the standard triangular decomposition of $\fg$, we  set $\mathfrak{m}=\mathfrak{t}\cap \mathfrak{n}^-$,    $\gb=\ga_N$, $F=F_\gb$ and $\gd =\ad_F$.
\\ \\
We remark that in \cite{M1} (3.8), the induction is set up in  an equivalent way using the Weyl group.  To explain this, let $s_i$ be the reflection corresponding to the simple root  $\ga_i$. If $w_i = s_i \ldots s_3s_2$ for $i\ge 2$, then $w_i \ga_1$ is the highest root of $\fk_i $. In particular if $\underline{v} = w_{N-1}$, 
 then $\gs = \underline{v} \ga_1$ and 
$\eta=w_N \ga_1$ are the highest roots of  $\fk$ and $\fg$ respectively.  Let $\gr$ be the half-sum of the positive roots of $\fg$, and define the dot action of $W$ on $\fh^*$ by $u\cdot \mu = u(\mu +\gr) - \gr$, for  $u \in W$ and $\mu\in\fh^*$.
Note that \be \label{120z} \eta = \gs + \gb,\ee
and $$\underline{v}^{-1}\eta = \eta,$$ so 
\be \label{120s} \underline{v}^{-1}\gb = \ga_{2} + \ldots +\ga_N.\ee
Now to ensure that $p$ in Lemma \ref{1768} is a positive integer (and likewise for previous steps in the induction), we need an assumption on $\mu$.  To ensure this assumption holds, we introduce the set 

\be \label{1tar}\Lambda =
\{\nu \in {\mathcal H}_{\ga_1,m}\mid \mbox{for } i \in [N], 
(\nu + \rho,\alpha_i^\vee) \in \mathbb{N} \backslash \{0\} \}
.\ee
Note that $\Lambda$ is Zariski dense in ${\mathcal H}_{\ga_1,m}$. 
If $\mu \in \underline{v}\cdot\Lambda$, then 
$$p= (\mu + \rho,\gb^\vee) = (\underline{v}^{-1}(\mu + \rho),(\underline{v}^{-1}\gb)^\vee)$$ is a positive integer by \eqref{120s}.
\bl \label{r25}  For $x\in U(\fm)$, there exists $n\ge 0$ such that $\gd^n x =0.$ 
\el
\bpf Let $N  = \{x\in U(\fm) \mid \gd^n x =0  \mbox{ for some } n >0\}.$ By \eqref{R24}, $N$ is a subalgebra of $U(\fm)$.  Thus it suffices to show that $N$ contains the generators  of $U(\fm)$. 
This follows from~\eqref{R7}.\epf
\noi 
\noi The next result follows by an easy computation.  Compare \cite{Jo1} 1.2.13 (1). 
\bl\label{lax} For $1\le i\le \ell -1$, we have 

$$v^{ i}\left[ \begin{array}{c}\ell -1\\i\end{array}\right]_{v}+v^{-(\ell -i)}\left[ \begin{array}{c}\ell -1\\i-1\end{array}\right]_{v} = \left[ \begin{array}{c}\ell \\i\end{array}\right]_{v}.$$ Equivalently, 
\be \label{2tar}
v^{ {-i(\ell -1-i)}}\left[ \begin{array}{c}\ell -1\\i\end{array}\right]_{v}+v^{ {-(i+1)(\ell -i)}}\left[ \begin{array}{c}\ell -1\\i-1\end{array}\right]_{v}  =v^{ {-i(\ell -i)}}\left[ \begin{array}{c}\ell \\i\end{array}\right]_{v}.\ee
 \el
\bl \label{3tar} For $k\ge 0$ and $z\in U$, we have $$\sigma (\mbox{\rm ad}_F^{k}( z))=v^{ {-2k}}\mbox{\rm ad}_F^{ k}(\sigma (z)).$$ 
 \el \bpf We can assume that  $z\in U^\eta$.  Then $\mbox{\rm ad}_F^{k}( z)\in U^{\eta-k\gb}$. Hence  the result follows from \eqref{cat1} and the fact that ${(\gb,-\gb)} = -2$.
\epf


\bl \label{Lx}
For $u\in U^{\nu }$, we have
\[
F^{\ell }u=\sum_{i=0}^{\ell}v^{ {-i(\ell -i)}}\left[ \begin{array}{c}\ell \\i\end{array}\right]_{v}v^{i(\beta ,\nu )}\mbox{\rm ad}_F^{\ell -i}(u)F^i.
\]
\el

\begin{proof}
Since $\sigma^i(u)=v^{i(\beta ,\nu )}u$, the given identity is the same as the identity
\begin{equation}\label{eq61}
F^{\ell }u=\sum_{i=0}^{\ell}v^{ {-i(\ell -i)}}\left[ \begin{array}{c}\ell \\i\end{array}\right]_{v}\mbox{\rm ad}_F^{\ell -i}(\sigma^i(u))F^i,
\end{equation}
and we prove this by induction on $\ell $. When $\ell =1$, this follows by writing Equation \eqref{R8} in the form
\be \label{eq6}Fu=\mbox{\rm ad}_F(u) +\sigma(u)F.\ee
 Suppose by induction that 
\be\label{eq7} y:= 
F^{\ell -1}u=\sum_{i=0}^{\ell -1}v^{ {-i(\ell -1-i)}}\left[ \begin{array}{c}\ell -1\\i\end{array}\right]_{v}\mbox{\rm ad}_F^{\ell -1-i}(\sigma^i(u))F^i.\ee
Then
{\small\[
\begin{aligned}
F^{\ell }u&=Fy=\sum_{i=0}^{\ell -1}v^{{-i(\ell -1-i)}}\left[ \begin{array}{c}\ell -1\\i\end{array}\right]_{v}F\cdot \mbox{\rm ad}_F^{\ell -1-i}(\sigma^i(u))F^i\\
&=\sum_{i=0}^{\ell -1}v^{ {-i(\ell -1-i)}}\left[ \begin{array}{c}\ell -1\\i\end{array}\right]_{v}\{\mbox{\rm ad}_F^{\ell -i}(\sigma^i(u))+\sigma (\mbox{\rm ad}_F^{\ell -1-i}(\sigma^i(u)))F\}F^i 
\\ 
&=\sum_{i=0}^{\ell -1}v^{ {-i(\ell -1-i)}}\left[ \begin{array}{c}\ell -1\\i\end{array}\right]_{v}\mbox{\rm ad}_F^{\ell -i}(\sigma^i(u))F^i+\sum_{i=0}^{\ell -1}v^{ {-(i+2)(\ell -1-i)}} \left[ \begin{array}{c}\ell -1\\i\end{array}\right]_{v}\mbox{\rm ad}_F^{\ell -1-i}(\sigma^{i+1}(u))F^{i+1}\\
&=\sum_{i=0}^{\ell -1}v^{ {-i(\ell -1-i)}}\left[ \begin{array}{c}\ell -1\\i\end{array}\right]_{v}\mbox{\rm ad}_F^{\ell -i}(\sigma^i(u))F^i+\sum_{i=1}^{\ell }v^{ {-(i+1)(\ell -i)}}\left[ \begin{array}{c}\ell -1\\i-1\end{array}\right]_{v}\mbox{\rm ad}_F^{\ell -i}(\sigma^{i}(u))F^{i}\\ 
&=\mbox{\rm ad}_F^{\ell }(u)+\sum_{i=1}^{\ell -1}\{v^{ {-i(\ell -1-i)}}\left[ \begin{array}{c}\ell -1\\i\end{array}\right]_{v}+v^{ {-(i+1)(\ell -i)}}\left[ \begin{array}{c}\ell -1\\i-1\end{array}\right]_{v}\}\mbox{\rm ad}_F^{\ell -i}(\sigma^i(u))F^i+\sigma^{\ell}(u)F^{\ell}\\
&= \sum_{i=0}^{\ell}v^{ {-i(\ell -i)}}\left[ \begin{array}{c}\ell \\i\end{array}\right]_{v}\mbox{\rm ad}_F^{\ell -i}(\sigma^i(u))F^i,
\end{aligned}
\]}where the third equality is due to \eqref{eq6} with $u$ replaced by $\mbox{\rm ad}_F^{\ell -1-i}(\sigma^i(u))$, the fourth 
follows from Lemma \ref{3tar},
the fifth follows a shift in the summation index for the second sum, the sixth is a rearrangement of terms,  and the last equality follows from \eqref{2tar}.
\end{proof}

\bl \label{lemX}
Use the same notation as Lemma~\ref{r25}. For $u\in U(\fm)^\nu$ there exists $k\in \mathbb{N}$ such that for all $\ell \in \mathbb{N}$,
\[
F^{\ell }u=\sum_{i=0}^{k}v^{ {-i(\ell -i)}}v^{(\ell -i)(\beta ,\nu )}\left[ \begin{array}{c}\ell \\i\end{array}\right]_{v}\mbox{\rm ad}_F^i(u)F^{\ell -i}.
\]
\el

\begin{proof}
By Lemma~\ref{r25} there exists $k\in \mathbb{N}$ such that $\delta^{k+1}(u)=0$.
The rest follows from Lemma \ref{Lx}, where the summation index $i$ is replaced by $\ell -i$.
\end{proof}
\bc \label{325} $S = \{F^n\mid n \ge 0\}$ is an Ore set in 
$U (\fn^-)$. \ec
\bpf It is enough to prove the Ore condition for $F$ and an element $u \in U (\fm )$. To do this use Lemma~\ref{lemX} to adapt the proof of Lemma 3.1 in  \cite{M1}. \epf
\subsection{The approach of  Shapovalov.} \label{dxd} 
We need a basic fact about Verma modules for $U_q(\fg)$.
\bt
\label{763} For $\gl, \mu \in P,$  every nonzero element of $\Hom_{U_q(\sfg)}(M(\mu),M(\lambda))$  is injective.
\et 
\bpf The key point is that $U_q(\fg)$  is a domain, \cite{JL}  Proposition 4.10.  Using this we can adapt the proof in \cite{M} Theorem 9.3.1.
The result is  implicit in \cite{Jo1} 3.4.9.
\epf  \noi 
Shapovalov elements in $U_q(\fg)$ may be constructed inductively using the next Lemma, compare \cite{H} Section 4.13 or  \cite{M} Lemma 9.4.3 for the non-quantum case. We  adapt the  notation to be consistent with our induction set-up.
\begin{lemma}\label{1768}
With the same notation  as Subsection \ref{dcd}, set $$\mu = s_\gb  \cdot  \gl, \;\;\; 
\gk = m\eta, \;\;\; \gk' = m\gs.$$  Assume that
\begin{itemize}
 \item[{(a)}]  $p = (\mu + \rho, \gb ^\vee) \in \mathbb{N}
\backslash \{0\}$.
\item[{(b)}]  $\theta' \in
U(\fm)^{-\gk'}$ is such that $\omega = \theta'v_\mu \in
M(\mu)$ is a highest weight vector.
\end{itemize}
Then there is a unique $\theta \in U(\mathfrak{n^-})^{-\gk}$ such that
\begin{equation} \label{21ndch76}
F^{p + m} \theta'v_\mu  = \theta F^p v_\mu.
\end{equation}
\begin{itemize}
\itemc 
$\theta F^pv_\mu $ is a highest weight vector in $U_q(\fg)F^pv_\mu  =
M(\gl)$.
\end{itemize}
\end{lemma}
\bpf  The proof in the classical case is based on the representation theory of $\fsl(2)$. Thus in the quantum case we use  the representation theory of $U_v(\fsl(2))$. Suppose that $E, F, K \in U_v(\mathfrak{g })$ satisfy the relations of  
$U_v(\fsl(2)),$ \cite{J} 1.1.  By \cite{Jo1}, Lemma 4.2.6, $v_\gl  =F^pv_\mu$ is a highest weight vector with weight $\gl$ in $M(\mu)$. Thus $U_q(\fg)F^pv_\mu  \cong 
M(\gl)$ by Theorem \ref{763}.  This accounts for the identification in (c). Set $\omega =\theta'v_\mu$.  
We claim that
\begin{equation} \label{21ndch75}
F^{p + m}\omega  \in  M(\gl).
\end{equation}
By Lemma \ref{lemX}, there is a positive integer $\ell$ such that
$F^{\ell}\omega  
 \in U(\mathfrak{n}^- )F^p$, and
hence $F^\ell \omega  \in M(\gl)$. We may assume
that $\ell \geq (p + m).$
Note that $\omega $ 
has weight $\mu - \gk'$ and
\be \label{11c} (\mu -\gk', \gb^\vee) = p +m -1.\ee 
The first equality below is \cite{J} 1.3 (5).
\by\label{drt} EF^\ell \omega &=& (F^\ell E  + [\ell]_vF^{\ell-1} [K;1-\ell])\omega \nn\\ &=&  [\ell]_v [1-\ell+(\mu - \gk',\gb )]_v F^{\ell-1} \omega \nn\\ &=&  [\ell]_v [p +m-\ell]_v F^{\ell-1} \omega .\ey
The second equality follows from \eqref{cat7} and the third from \eqref{11c}. 
 If $\ell> p + m$, then the coefficient of  $F^{\ell-1}\omega $  in \eqref{drt} is a non-zero Laurent polynominal in~$v$.  Thus ~${M}(\gl)$ contains $F^{\ell-1}\omega $.  Repeating this argument gives (\ref{21ndch75}).
Since~$M(\gl)
=U(\mathfrak{n}^-)e^pv_\gl$ is a free $U(\mathfrak{n}^-)$-module,
there exists an element $\theta \in U(\mathfrak{n}^-)^{- \gk }$ such
that (\ref{21ndch76}) holds. Since $F$ is not a zero
divisor in $U(\mathfrak{n}^-)$, the element $\theta$ satisfying~(\ref{21ndch76}) is unique.
\\ \\
To show (c), write 
$v_\gl=F^pv_\mu$ By Theorem \ref{763}, the submodule of $M(\mu)$ generated by $v_\gl$ is isomorphic to $M(\gl)$. 
We must show that $E_\ga v_\gl=0$ for every simple root $\ga$.  If $\ga\neq \gb$ this follows since $[E_\ga,F_\gb] =0$. If $\ga= \gb$ 
we have by \eqref{21ndch76}, 
$F^{p + m}\omega = \theta v_\gl$ 
 and taking $\ell = p +m$ in \eqref{drt}, this is zero because $[p +m-\ell]_v =0$.\epf

\subsection{Comparison of two approaches to Shapovalov elements.} \label{2cal}
Recall that the subdiagonal entries in the matrix \eqref{coc} are $-c_i$, where $c_i =h_i(\gl)$, see  
\eqref{0q}. 
\bl\label{7cal}
Assume  $(\gl +\gr, \eta) = 1$ and $(\gl +\gr, \gb) = -p$ as in Lemma \ref{1768}. Then 
\bi \itema $$c_{N-1} = h_{N-1}(\gl) = -q^{-1}v^{-p}[p+1]_v.$$
\itemb For  $ i \in  [{N-2}]$, we have $c_i= h_i(\mu)$.
\ei\el
\bpf 
From \eqref{120z} and the hypotheses,
$(\gl +\gr, \gs) = p+1$.  Therefore using \eqref{1q} and \eqref{0q} we obtain (a).\\
\\
If $ i \in  [{N-2}]$,  $c_i$ in 
 depends only on the value of $(\gl +\gr,\gs_i)$. Since $(\gb,\gs_i)=0,$ we have $(\mu +\gr,\gs_i)=(\gl +\gr,\gs_i)$. Hence 
$h_i(\mu)=h_i(\gl)$.
\epf
We need an analog of \eqref{120} for $U_q(\fsl(N))$.  Thus define
\be \label{120v} \J = \{J\subseteq [N]\mid 1, N\in J\}. \ee  
For $J\in\J$, set $f_J$ as in \eqref{ddd}. If $J$ is as in 
\eqref{dod}, with   
$j_0 =  1$ and $j_{s+1} =N$, define 
$J_1, J_2\in\I$ by
$$J_1= \{j_0 , j_1, \ldots ,j_s, N, N+1\} \mbox{ and } J_2=     \{j_0 , j_1, \ldots ,j_s,  {N+1} \}.$$ 
Note that by \eqref{cat0},
\be \label{120w} \deg f_J = \gep_N - \gep_1. \ee
By \eqref{dog},
 \eqref{cat1} and \eqref{R8},  with $F=f_{N,N+1}$, 
\by\label{dgt} \ad_{F }(f_{i,N}) &=& f_{N,N+1} f_{i,N}  -v f_{i,N} f_{N,N+1} = -qf_{i,N+1}.\ey
More generally, we have the part (a) of the following Lemma. Part (b) is a similar relation.
\bl \label{dkt} 
\begin{enumerate}
\itema
 $\ad_{F }(f_J) =  -qf_{J_2}.$
\itemb $f_JF = f_{J_1}$.
\end{enumerate}
\el
\bpf
(a) Induct on $s=|J|$. The case $s=1$ follows from \eqref{dgt}. Suppose $s>1$ and write $f_J=f_{j_0,j_1}f_{J'}$ and $f_{J_2}=f_{j_0,j_1}f_{J_2'}$.
By induction, $\ad_{F }(f_{J'}) =  -qf_{J_2'}$. Left multiplication of it by~$f_{j_0,j_1}$ clearly is $f_{j_0,j_1}\ad_{F }(f_{J'}) =  -qf_{J_2}$. The left side here is $f_{j_0,j_1}\ad_{F }(f_{J'})=f_{j_0,j_1}(Ff_{J'}-Kf_{J'}K^{-1}F)$, where $K=K_{N,N+1}$. Since~$\{j_0,j_1\}$ and~$\{N,N+1\}$ are disjoint, $f_{j_0,j_1}$ commutes with both $F$ and $K$, and the last expression becomes $Ff_{j_0,j_1}f_{J'}-Kf_{j_0,j_1}f_{J'}K^{-1}F=Ff_{J}-Kf_{J}K^{-1}F=\ad_{F }(f_J)$, as needed. 
Part (b) is proved in a similar way. The case  $s=1$ of the induction holds by the definition of $f_{J_1}$.
\epf
\noi If $u= f_J$, we can take $k= 1$ in Lemma 
 \ref{lemX}. Note that $(\gb, \deg f_J) = 1$ by \eqref{120w}. Hence by Lemmas \ref{dkt} and \ref{7cal}, 
\by\label{diit}
F^{p+1 }f_J&=&\sum_{i=0}^{1}v^{ {-i(p+1 -i)}}v^{(p+1 -i)(\beta ,\nu )}\left[ \begin{array}{c}p+1 \\i\end{array}\right]_{v}\mbox{\rm ad}_F^i(f_J)F^{p+1 -i}.
\nn\\
&=&   v^{p+1}f_J F  +[p+1]_v \ad_{F }(f_J))
F^{p}\nn\\
&=&  v^{p+1}( f_{J_1} - q^{-1}v^{-p}[p+1]_vf_{J_2})
F^{p}\nn\\
&=&  v^{p+1}( f_{J_1} +  c_{N-1}f_{J_2})
F^{p}.\ey
The following result relates the two approaches to Shapovalov elements. 

\bt \label{doot} $$F^{p+1}{\stackrel{\longrightarrow}{{\rm det}}}\;\cD^{N}(\mu) = v^{p+1}
{\stackrel{\longrightarrow}{{\rm det}}}\;\cD^{N+1}(\gl)F^p.$$ 
\et
\bpf Let $D_2$ be the cofactor of entry $-c_{N-1}$ in $\cD^{N+1}(\gl)$, see \eqref{coc}, and let $D_1$ be the cofactor of entry $f_{N,N+1}$.  
It follows from Lemma \ref{7cal}, that 
\be \label{120a}\sum_{J\in\J} H_J(\mu)f_J= {\stackrel{\longrightarrow}{{\rm det}}}\; D_1 = {\stackrel{\longrightarrow}{{\rm det}}}\;\cD^{N}(\mu).\ee 
Now  $f_{J_2}$ is obtained from $f_{J}$ 
by replacing the last factor $f_{j_s,N}$ 
by $f_{j_s,N+1}.$ Since $D_2$ is obtained from $D_1$ by replacing the entries $f_{i,N}$ in the last column by $f_{i,N+1}$, it follows that  
\be \label{120b}\sum_{J\in\J} H_J(\mu)f_{J_2}= {\stackrel{\longrightarrow}{{\rm det}}}\; D_2 
.\ee 
Thus from \eqref{diit}, \eqref{120a} and \eqref{120b},
\by F^{p+1}{\stackrel{\longrightarrow}{{\rm det}}}\;\cD^{N}(\mu) &=& F^{p+1}\sum_{J\in\J} H_J(\mu)f_{J}\nn\\ 
&=& v^{p+1}\sum_{J\in\J} H_J(\mu)(f_{J}F +  c_{N-1}f_{J_2})
F^{p}
\nn\\
&=&v^{p+1}( {\stackrel{\longrightarrow}{{\rm det}}}\;D_1 f_{N,N+1}+c_{N-1}{\stackrel{\longrightarrow}{{\rm det}}}\;D_2)F^p\nn\\ &=&v^{p+1}{\stackrel{\longrightarrow}{{\rm det}}}\;\cD^{N+1}(\gl)F^p,\nn 
\ey where the last equality uses cofactor expansion along the last row of $\cD^{N+1}(\gl)$.
\epf \noi
 For a Shapovalov element, there is a normalization condition that a certain coefficient in $\H$ is equal to one.  No such condition is made on the elements $\theta'$ and $\theta$ in Lemma \ref{1768}.
From   \eqref{diit}, we can easily see what these coefficients are in Theorem  \ref{doot}. 
In \eqref{120v}, take $J = [N]\in \J$. Now $f_J = f_{[N]}$ and $f_{J_1}= f_{[N+1]}$ occur with coeffiecient 1, in the Shapovalov elements $\theta_\gs$ and  $\theta_\eta$ respectively and we have,   
\be\label{daat}
F^{p+1 }f_J = v^{p+1}f_{J_1}F^{p}.\ee 
This means that if we take $m=1$ in  \eqref{21ndch76} and define $\overline{\theta} = v^{-p-1}\theta$,  we have 
\begin{equation} \label{23ndch76}
F^{p + 1} \theta'v_\mu  = 
v^{p+1}\overline{\theta} 
 F^p v_\mu = v^{p+1}\overline{\theta} 
  v_\gl.
\end{equation}
This means that if $\theta'(\mu)$ is the evaluation of the Shapovalov element 
for $\gs$, then    
$\overline{\theta}(\gl)$ is the evaluation of the Shapovalov element 
for $\eta$, in accordance with Theorems \ref{T} and \ref{doot}.

\subsection{A uniform construction for Shapovalov elements in $U(\fb^-)$.} \label{nir} 

Up to this point we have only evaluated a Shapovalov element
$\theta_{\gamma,m}$ at points $\gl \in \cH_{m,\gc}$.
However, to study the behavior of powers of  Shapovalov elements,   we need to evaluate at arbitrary points, $\gl \in \fh^*$.
However some care must be taken since the Shapovalov element $\theta_{\gamma,m}$ is only defined modulo the ideal
$U({\mathfrak b}^{-} ){\mathcal I}(\mathcal H_{\gc, m})$.
\\ \\
To explain the problem 
we briefly review the construction of  Shapovalov elements from \cite{M} 9.4, adapted to the quantum case.  
First write the root $\gc$  in the form $\gc =w \gb$ with $\gb$ a simple root and $w\in W$.  Next  by induction on the length of  the $w$, we construct 
$\theta_{\gamma,m} \in U_q(\fn^-)^{-m\gc}$ such that $\theta_{\gamma,m} v_\gl$ is a highest weight vector for $M(\gl)$ whenever $\gl\in {\mathcal H}_{\gc,m}$,  
see the remarks after \eqref{boo}. 
Then we fix a (non-unique) lifting of $\theta_{\gamma,m}$  to an element $\Theta_{\gamma,m}\in U_q(\fb^-)^{-m\gc}= U_q(\fn^-)^{-m\gc} \ot \H$ such that $
\Theta_{\gamma,m} v_\gl= \Theta_{\gamma,m}(\gl) v_\gl =\theta_{\gamma,m} v_\gl$,
for all $\gl \in {\mathcal H}_{\gc,m}$. 
\\ \\
\noi  The problem is that  if $\gl \in \mathcal H_{\gc, m}$, then $\Theta_{\gamma,m}v_\gl$ is a highest weight vector of weight  $\gl-m\gc$,  but $\gl-m\gc\notin \mathcal H_{\gc, m}$, and thus we need 
to evaluate $\Theta_{\gamma,m}$   at  points that are not in 
$\mathcal H_{\gc, m}$. 
Thus the non-uniqueness of the lifting of $\theta_{\gamma,m}$ presents a potential problem. To resolve this issue, we give a uniform inductive  construction 
of $\Theta_{\gamma,m}$, following \cite{M1} 5.1.2 or \cite{M3}. This depends on a specific choice of $\gb$ and a shortest length expression for $w$.
\\ \\
With the notation  of Lemma \ref{r25}, let  $S$ be the multiplicative subset $S=\{F^n|n\ge 0\}$ as in Corollary 
\ref{325}, 
and let $U_S$ be the corresponding Ore localization. 
We define generalized conjugation operators for non-isotropic roots. 
Lemma \ref{lemX} leads to a formula for conjugation by $F^r$ which extends to the
$1$-parameter family of automorphism $\Psi_r, r\in \ttk$ 
of $U_S,$ given by
\be \label{cows} \Psi_{r}(u) = 
\sum_{i
 \ge 0} v^{ {-i(r -i)}}\left[ \begin{array}{c}r \\i\end{array}\right]_{v}v^{(r - i)(\beta ,\nu )}\mbox{\rm ad}_F^{i}(u)F^{i -r},\ee
for $u\in U^{\nu }$. 
Note that for $r\in\N$  and $u\in U$  we have $\Psi_r(u) = F^r u F^{-r}$ by 
Lemma~\ref{Lx}.
We extend $\Psi_r$ to an automorphism of 
$U=U (\fn_{i} )$ by setting 
  $\Psi_r(F) = {F}$. Then,  
\be\label{yam} \Psi_r(FuF^{-1}) = \Psi_{r+1}(u). \ee
As noted earlier, the value of $\Psi_r(u)$ for $r\in \ttk$   is determined by its values for $r\in\N$. In particular the two previous equations hold for all $r\in \ttk.$ To stress the dependence of $\Psi_r$ on the root $\gb$ we sometimes denote it by $\Psi^\gb_r.$
\\ \\   Because of the polynomial nature of the coefficients in 
\eqref{cows}, this implies that $\Psi_{r}$ is an automorphism of $U_S$ for all $r\in \ttk.$  Now 
\eqref{21ndch76} 
implies that $\Psi_{r}(e^{m}_{- \alpha}\gth_{\gamma',m}) = \gth_{\gamma,m}$ if 
$\mu \in \cH_{\gc',m}$ and $r = (\mu + \rho, \alpha^\vee)$. 
 
\bl
Suppose that  $\mu \in \cH_{\gc',m}, \;\ga$ is a simple  root, $\gc = s_\alpha\gc',\;\;  (\gc, \alpha^\vee)=1$ and $\gl = s_\alpha\cdot \mu.$
Assume  $r = (\mu + \rho, \alpha^\vee)$. 
Then  we have
\be \label{bam1} 
\theta_{\gamma,m}(\gl)= \Psi_{r}(F^{m}_{- \alpha}\theta_{\gamma',m}(\mu)).\ee\el
\bpf If $r$ is a positive integer, this follows from Lemma \ref{1768}. Since both sides depend polynomially on $\mu$ we have the result.  \epf


\noi 
We propose \eqref{bam1} as the definition of $\theta_{\gamma,m}(\gl)$ for  arbitrary points, $\gl = s_\alpha\cdot \mu \in \fh^*$, without the assumption that $\mu \in \cH_{\gc',m}.$ Thus $\theta_{\gamma,m}$ is now a function from $\fh^*$ to $U(\fn^-)^{-m\gc}$  which evaluates correctly on $\gl \in {\mathcal H}_{\gc,m}$ and has $H_{\pi^0}=1$  as the coefficient of 
$e_{-\pi^0}$.
Thus the Shapovalov elements  $\Gt_{\gamma,m}\in U(\fb)^{-m\gc}$, are defined inductively 
by 
\be \label{bam}\Gt_{\gamma,m}(\gl)= \Psi_{r}(e^{mq}_{- \alpha}\Gt_{\gamma',m}(\mu)),\ee
 where $r = (\mu + \rho, \alpha^\vee)$. Now we have a uniform construction  
of such elements. In the case when $\gb$ is a simple root, we set $\Gt_{\gb, m} = F_{-\gb}^m\in U(\fb)^{-m\gb}$.
\\ \\
At first glance it appears that we have only 
$\Gt_{\gamma,m}(\gl)\in U_S$, where $S$ is the Ore set from Lemma \ref{325}. However we have
\bl \label{c4} For all $\gl \in \fh^*,$ and $j>0,$ we have $  \Gt_{\gamma,j}(\gl)\in U$.\el%
\bpf It is easy to adapt the proof from  
\cite{M1} Corollary 5.6.
\epf

\subsection{Powers of Shapovalov elements.} \label{1zzprod}
\bt \label{1calu} If $\gl \in \cH_{\gc,m}$, then
\be\label{1CL}
 \Gt_{\gamma,m}(\gl) =\Gt_{\gamma,1}(\gl-(m-1)\gc)\ldots \Gt_{\gamma,1}(\gl-\gc)\Gt_{\gamma,1}(\gl).
\ee
Equivalently if $v_\gl$  is a highest weight vector of weight $\gl$ in the Verma module $M(\gl)$ and $\gl \in \cH_{\gc,m}$, we have $\Gt_{\gamma,m}v_\gl = \Gt_{\gamma,1}^m v_\gl$. If $\gl \in \cH_{\gc,m}$ this is a highest weight vector in $M(\gl)$.
\et
\bpf 
Clearly \eqref{1CL} holds if $\gc$ is a simple root. Suppose that 
\be \label{121c}\gl = s_\alpha\cdot \mu,\;\;\gc = s_\alpha\gc',\;\;p = (\mu + \rho, \alpha^\vee),\;\;q = (\gc, \alpha^\vee).\ee
For $i=0,\ldots, m-1$ we have $$(\mu + \rho-i\gc', \alpha^\vee)= p+i,$$ so by the inductive definition \eqref{yam} and \eqref{bam},

\by \label{1inp}
\Gt_{\gamma,1}(\gl-i\gc)&=& \Psi_{ p + i}(F\Gt_{\gamma',1}(\mu-i\gc'))\nn\\
&=& \Psi_{ p }(F^{(i+1)}\Gt_{\gamma',1}(\mu-i\gc')F^{-i}).\ey
Now using the corresponding result for $\Gt_{\gamma',m}(\mu)$ we have
\by  \label{tea}F^{m} \Gt_{\gamma',m}(\mu) &=& F^{m}\Gt_{\gc',1}(\mu-(m-1)\gc')\ldots \Gt_{\gc',1}(\mu-\gc')\Gt_{\gc',1}(\mu)\nn\\
 &=&
F^{(m-1)}(F\Gt_{\gc',1}(\mu-(m-1)\gc'))F^{-(m-1)}\times\nn\\
&&\times F^{(m-2)}(F^{}\Gt_{\gc',1}(\mu-(m-2)\gc'))F^{-(m-2)}\times\nn\\
&\cdots & \nn\\
& \cdots &  \times F^{(i+1)}\Gt_{\gamma',1}(\mu-i\gc')F^{-i}\times\nn\\
&\cdots & \nn\\
& &\times \;F^{}(F^{}\Gt_{\gamma',1}(\mu-\gc'))F^{-1}\times F^{}\Gt_{\gamma',1}(\mu).\nn\ey
The result follows by applying the automorphism $\Psi_p$ to both sides and using \eqref{bam} and \eqref{1inp}.\epf

\bibliographystyle{plain}

\end{document}